\documentclass[a4paper]{article}

\usepackage{amsfonts}
\usepackage{amsopn}
\usepackage{amsmath}
\usepackage{amsthm}
\usepackage{latexsym}
\usepackage{amssymb}
\usepackage{dsfont}
\usepackage{graphicx}

\newcommand{\EE}{{\mathbb{E}}}

\newcommand{\PP}{{\mathbb{P}}}
\renewcommand{\P}{{\PP}}
\newcommand{\RR}{{\mathbb{R}}}

\newcommand{\nF}{{\mathcal{F}}}

\newcommand{\re}{{\mathrm{e}}}

\newcommand{\bP}{{\boldsymbol{P}}}

\newcommand{\p}{\partial}
\newcommand\1{\mathds{1}}

\renewcommand{\phi}{\varphi}

\newcommand{\abs}[1]{\left|#1\right|}
\newcommand{\menge}[1]{\left\lbrace #1\right\rbrace }

\newcommand{\N}{\mathbb{N}}
\newcommand{\Z}{\mathbb{Z}}
\newcommand{\R}{\mathbb{R}}

\newcommand{\ten}{\fontsize{10pt}{\baselineskip}\selectfont}

\title{Transport in Stochastic Goupillaud Media}

\author{Michael Oberguggenberger\thanks{Unit of Engineering Mathematics, University of Innsbruck,
Technikerstra\ss e 13, 6020 Innsbruck,
Austria, (michael.oberguggenberger@uibk.ac.at)}
 \and
Martin Schwarz\thanks{Unit of Engineering Mathematics, University of Innsbruck,
Technikerstra\ss e 13, 6020 Innsbruck,
Austria, (mr2cef@gmail.com)}
}
\date{}

\begin{document}
\maketitle

\begin{abstract}
  The paper addresses one-dimensional transport in a Goupillaud medium (a layered medium in which the layer thickness is proportional to the propagation speed), as a prototypical case of wave propagation in random media. Suitable stochastic assumptions and limiting procedures lead to characteristic curves that are L\'evy processes. Solutions corresponding to the discretely layered medium are shown to converge to limits as the thickness of the layers goes to zero. The probability distribution of the limiting characteristic curves is explicitly computed and exemplified when the underlying L\'evy process is an inverse Gaussian process.
\end{abstract}

\section{Introduction}
\label{sec:intro}

The study of wave propagation in random media has a long history and a wealth of applications in the material sciences, geotechnics, and seismology (Fouque et al. 2007). It is common practice to model the stochastic properties of such media by means of random fields (Ghanem \& Spanos 1991). The present paper focuses on singular stochastic limits of one-dimensional layered media.
A prototypical case is given by Goupillaud media, which are characterized by the property that the travel time of a wave through each layer is constant. This type of media has been introduced by (Goupillaud 1961) to study seismic wave propagation. The theme has been pursued over the past decades, mainly focusing on describing transmission and reflection of waves (Velo \& Gazonas 2019). Random Goupillaud media generated from Markov chains as well as their limiting behavior as the layer thickness goes to zero have also been studied (Burridge et al. 1988). It seems to be less known that in unidirectional wave propagation (transport equations), stochastic Goupillaud media can be constructed possessing an extremely rough stochastic structure, of a spatially much higher irregularity than e.g. in (Flandoli 2011).

This paper will clarify the case of one-dimensional transport under assumptions that will lead to characteristic curves given by an increasing L\'evy process.

One-dimensional transport is described by the equation
\begin{equation} \label{eq:trans}
\p_t u(t,x) + c(x) \p_x u(t,x)= 0
\end{equation}
with initial data $u(0,x) = u_0(x)$. If the transport velocity c(x) and the initial data are sufficiently regular, the solution is given by
\begin{equation} \label{eq:solu}
u(t,x) = u_0(\gamma(x,t;0))
\end{equation}
where $\gamma(x,t;0)$ is the point of intersection of the characteristic curve through $(x,t)$ with the $x$-axis. The characteristic curves are obtained as solutions to the ordinary differential equation (derivative with respect to $\tau$)
\begin{equation*} 
   \dot{\gamma}(x,t;\tau) = c(\gamma(x,t;\tau)),\quad \gamma(x,t;t)=x.
\end{equation*}
The material properties of the medium are encoded in the transport speed $c(x)$.
The Goupillaud assumption is that $c(x)$ has the constant value $c_k$ in the $k$th layer; the thickness and the travel time are related through $\Delta x_k = c_k\Delta t$.

Further, the propagation speeds $c_k$ (and hence the $\Delta x_k$) will be given by independent, identically distributed random variables. At this stage, various choices of the type of random variables as well as scalings are possible. For the wave equation, such scalings leading to fairly regular limiting processes have been introduced in (Burridge et al. 1988) and studied in (Fouque et al. 2007, Nair \& White 1991).
Setting $\Delta t = 2^{-N}$, it is assumed here that the $(N+1)$st layering is a dyadic refinement of the $N$th layering and that, at each stage $N$, the $c_k^{(N)}$ (and hence the $\Delta x_k^{(N)}$) are positive, independent and identically distributed random variables.

The procedure of dyadic refinements on the time axis leads to infinitely divisible, positive random variables, which can be constructed as increments of a strictly increasing L\'evy process, which in turn defines the limiting characteristic curve of Equation\;\ref{eq:trans} through the origin as $N\to\infty$.

The following results will be presented: The characteristic curves converge to translates of the path of the said L\'{e}vy process with probability one and at almost all $x$ and $t$. In addition, the piecewise solutions $u^{(N)}(t,x)$, obtained by solving the transport equation for the Goupillaud medium at stage $N$, converge to a limiting stochastic process almost surely and in the $p$th mean; the distribution of the limiting characteristic curves can be explicitly computed. The calculations will be made explicit when the underlying L\'evy process is an inverse Gaussian process. The slightly more complicated calculations for a Poisson process with positive drift can be found in (Schwarz 2019).

The limiting solution process $u(t,x)$ is constant along the limiting characteristic curves, as in the case of classical transport. However, the limiting characteristics may possibly have infinitely many jumps on each interval. Due to this high degree of irregularity, one cannot give a meaning to the limiting process $u(t,x)$ as a solution to Equation\;\ref{eq:trans}; it is just a limit of piecewise classical solutions. This situation is quite common in the theory of singular stochastic partial differential equations, see e.g. (Hairer 2014).

The plan of the paper is as follows: In Section 2, the stochastic Goupillaud medium is set up and analyzed. Section 3 is devoted to the limiting behavior of the characteristic curves and the solutions. In Section 4, it is shown how the probability distribution of the limiting characteristic curves can be computed, and explicit formulas in the inverse Gaussian case are given. The final section contains the conclusions.

Detailed proofs of the results of Sections 2 and 3 have been presented in (Baumgartner et al. 2017). Section 4 is based on the PhD thesis of the second author (Schwarz 2019), where also the case of a Poisson process can be found.

Notation: In the sequel, $\Z$, $\N$, $\R$ denote the sets of integers, nonnegative integers, and real numbers, respectively.

\section{The Goupillaud medium}
\label{sec:Goupillaud}

If the transport velocity $c$ in Equation\;\ref{eq:trans} is constant the characteristic curves are simply given by $\gamma(x,t;\tau)=x+c(\tau-t)$. In the case of a layered medium, $c(x)$ is piecewise constant, and the characteristic curves are polygons. Assuming continuity across interfaces, the solution (Formula\;\ref{eq:solu}) is given as a continuous, piecewise differentiable function, which solves Equation\;\ref{eq:trans} in the weak sense.

\subsection{Dyadic deterministic structure}
\label{subsec:dyadic}

The discrete, deterministic Goupillaud medium is set up as follows. Choose a time step $\Delta t$, producing a sequence $t_j = j\Delta t$, $j = 0,1,2,\ldots$ of points in time. The layers are defined through a strictly increasing sequence $(x_k)_{k\in\Z}$ with $x_0=0$ and $x_k\to\pm \infty$ as $k\to \pm \infty$. Let $\Delta x_k = x_{k}-x_{k-1}$. The coefficient $c(x)$ is obtained as
\begin{equation*} 
c(x) = \Delta x_k/\Delta t \mbox{\ for\ }  x_{k-1} \leq x < x_{k},
\end{equation*}
and $k$ running from $-\infty$ to $\infty$. In other words, the time for passing a layer $\Delta x_k$ is constant, namely $\Delta t$.
Call $c_k$ the value of $c(x)$ in the $k$th layer, that is, for $x_{k-1}\leq x < x_{k}$. Then the Goupillaud relation
$\Delta x_k = c_k\Delta t$
holds for all $k$, with constant $\Delta t$. The structure of the Goupillaud medium makes computing the values of the characteristic curves $\gamma(x,t;\tau)$ in the grid points very simple. In fact,
\begin{align}\label{eq:grid0}
\gamma(x_k,t_l;t_j)  = x_{j+k-l}
\end{align}
for all integers $j,k,l$. Since every point $(t,x)$ is just a convex combination of the neighboring grid points, the values $\gamma(x,t;\tau)$ can be easily obtained anywhere.

\begin{figure}[htb] 	
\centering
\includegraphics[width = 0.9\linewidth]{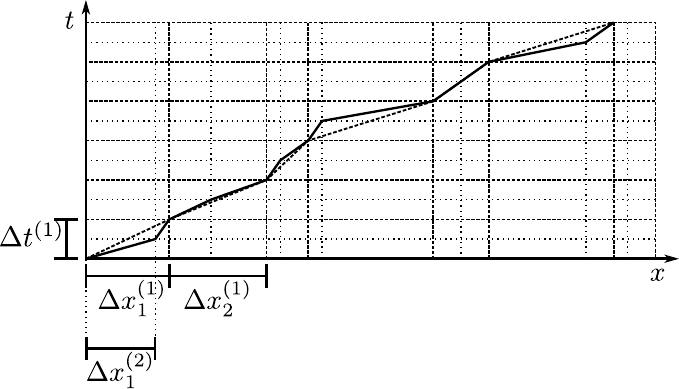}
	\caption{Figure 1. Illustration of Goupillaud medium and dyadic refinement. Shown are some layers of the random initial grid with constant time steps $\Delta t^{(1)}$ (dashed horizontal and vertical lines) and initial broken characteristic curve $x = \gamma^{(1)}(0,0,t)$ starting at the origin (dashed). Replacing $\Delta t^{(1)}$ by $\Delta t^{(1)}/2$ results in the second layering (dotted horizontal and vertical lines) and the characteristic curve $x = \gamma^{(2)}(0,0,t)$ (solid broken line), which passes through the old and new grid points.}
\label{fig:Gitter_fein}
\end{figure}

The next step is to set up a dyadic refinement of the initial grid. Define
\[
 \Delta t^{(N)}= 2^{-N},\quad t^{(N)}_j=j2^{-N}
\]
and let $x^{(N)}_k\in \R$, $k\in\Z$, be a strictly increasing sequence of spatial points (or equivalently, a sequence of propagation speeds $c^{(N)}_k > 0$ satisfying
$\Delta x^{(N)}_k = c^{(N)}_k \Delta t^{(N)}$). The requirement that each resulting grid is a dyadic refinement of the previous one means that
\begin{equation*} 
 \big(t^{(N+1)}_{2j},\, x_{2k}^{(N+1)}\big)=\big(t^{(N)}_{j},\, x_{k}^{(N)}\big).
\end{equation*}
This condition implies
\begin{equation*}
 \Delta x^{(N)}_k= \Delta x^{(N+1)}_{2k-1}+ \Delta x^{(N+1)}_{2k}.
\end{equation*}
A random realization of the initial grid and the first step of the dyadic refinement can be seen in Figure \ref{fig:Gitter_fein}, together with the corresponding broken characteristic curves starting at $(0,0)$.

Inductively, one obtains
\begin{align} \label{eq:consistency}
\Delta x_{k}^{(N)}= \sum_{i=1}^{2^M} \Delta x_{(k-1)2^M+i}^{(N+M)}
\end{align}
for all $N,M\in\N,k\in\Z$.
The value of the characteristic curve $\gamma^{(N)}$ in the grid points is readily obtained according to (\ref{eq:grid0}). For whatever integers $N,j,k,l$ it holds that
\begin{align}\label{eq:grid}
\gamma^{(N)}\left(x_k^{(N)},t^{(N)}_l;t_j^{(N)}\right)  = x^{(N)}_{j+k-l}.
\end{align}
For any  $N\in \N$ and $\tau\in [t_{k-1}^{(N)}, t_{k}^{(N)})$, the characteristic curve through the origin $\xi^{(N)}(\tau)=\gamma^{(N)}(0,0;\tau)$ can be represented as
\begin{equation}\label{eq:lambda}
\xi^{(N)}(\tau) = \alpha^{(N)} (\tau) \xi^{(N)}(t_{k-1}^{(N)})
                + \big(1-\alpha^{(N)} (\tau) \big)\xi^{(N)}(t_{k}^{(N)}),
\end{equation}
where $\alpha^{(N)}(\tau) =\big(t_{k}^{(N)}-\tau\big)2^N$ and $\xi^{(N)}(t^{(N)}_k)=x^{(N)}_k$ by Formula\;\ref{eq:grid}.

In other words, $\xi^{(N)}$ is an increasing polygon through $(t_k^{(N)},x_k^{(N)})$, $k\in\Z$.
For $(t,x)\in\R^2$ one obtains the characteristic curve through $(t,x)$ by
\begin{equation*} 
\gamma^{(N)}(x,t;\tau)= \xi^{(N)}\big(\tau+(\xi^{(N)})^{-1}(x)-t\big),
\end{equation*}
i.e., by shifting $\xi^{(N)}$ in time direction such that it passes through $(x,t)$.

\subsection{The stochastic model}
\label{subsec:stochmod}

This subsection is devoted to formulating the stochastic assumptions underlying the proposed model of a randomly layered medium, in which $\Delta x_k^{(N)}$ is random.
The decisive assumption is that for each $N\in\N$, the increments are positive, independent and identically distributed random variables $\Delta x_k^{(N)}$, $k\in\Z$.

Together with the previous consistency assumption (Formula\;\ref{eq:consistency}), this implies that $\Delta x_{k}^{(N)}$ is infinitely divisible for every $k\in\Z$ and $N\in\N$.
Let $\mu =\P\circ (\Delta x^{(0)}_1)^{-1}$ be the distribution of $\Delta x^{(0)}_1$. It follows from standard probabilistic arguments (Sato 1999, Section 7) that $\Delta x_{k}^{(N)}\sim \mu^{*1/2^N}$, $k\in\N$, the $2^N$th unique root of $\mu$, and that there exists a L\'evy process $L=(L(t))_{t\in\R}$ on a probability space $(\Omega,\nF,\P)$, with $\P\circ L(1)^{-1}=\mu$, that is, $\Delta x^{(0)}_1$ has the same distribution as $L(1)$.
These conditions are met by Poisson or Gamma processes with positive drift, as well as by the inverse Gaussian process, for example.

In other words, the condition of independent and identically distributed random variables, together with the dyadic structure, leads to a limiting L\'evy process. Conversely, given a L\'evy process, its dyadic increments produce sequences of independent identically distributed random variables. Thus one may equivalently use a given L\'evy process as starting point for defining the stochastic Goupillaud medium.

Indeed, take a L\'evy process $L$, let $t_k^{(N)} = k2^{-N}$ be as in Subsection~\ref{subsec:dyadic} and define
\begin{align*}
&\Delta x_{k}^{(N)}(\omega)= L\big(\omega,t_k^{(N)}\big)-L\big(\omega,t_{k-1}^{(N)}\big),\\
   &x_{k}^{(N)}(\omega) = \sum_{i=1}^k \Delta x_{i}^{(N)}(\omega) = L\big(\omega,t_k^{(N)}\big)
\end{align*}
for $k>0$ and similarly for $k\leq 0$. The consistency condition (Formula\;\ref{eq:consistency}) is clearly satisfied.
Let furthermore $L^{(N)}(\omega,\cdot)$ be the piecewise affine interpolation of $L(\omega,\cdot)$ through the grid points $(t_k^{(N)},x_k^{(N)}(\omega))$ as in Equation\;\ref{eq:lambda}. This construction is carried out pathwise for fixed $\omega\in\Omega$.

\section{Limits as the time step goes to zero}
\label{sec:convergence}

The main result of this section is that the characteristic curves of the discrete Goupillaud medium converge to limiting curves (almost surely almost everywhere). This will imply that the solutions to the transport equation converge to a limit as well (in a sense to be made precise). The crucial observation is that the paths of a L\'evy process are c\`adl\`ag (continue \`a droite, limite \`a gauche) almost surely, i.e., they are continuous from the right and have left-hand limits.

\subsection{Convergence of characteristic curves}
\label{subsec:convchar}

Fix an increasing L\'evy process $L$ as considered in Subsection\;\ref{subsec:stochmod} and let $L^{(N)}(\omega,t)$ be the piecewise affine interpolations of its paths.
At level $N$, the L\'evy process defines a discrete Goupillaud medium. The corresponding characteristic curves are given by inserting $L$ in place of $\xi$ in the formulas at the end of Subsection\;\ref{subsec:dyadic}, namely
$
\gamma^{(N)}(\omega; x,t;\tau)=L^{(N)}\big(\omega;\tau+(L^{(N)}(\omega))^{-1}(x)-t\big).
$
It will be shown that these curves converge, as $N\to\infty$, to
\begin{equation}\label{eq:charcurves}
\gamma(\omega;x,t;\tau)=L\big(\omega;\tau+L(\omega)^*(x)-t\big)
\end{equation}
where $L(\omega)^*$ is defined through
\begin{equation}\label{eq:Lstar}
L(\omega)^*(x)= \inf\menge{t\in\R: L(\omega,t)\geq x}.
\end{equation}
For almost all $\omega\in\Omega$, the map $t\to L(\omega,t)$ is an increasing c\`adl\`ag function. Using this property, one can derive that, at fixed $(x,t)\in \R^2$,
\[
\lim_{N\to\infty}\gamma^{(N)}(\omega;x,t;\tau_0)= \gamma(\omega;x,t;\tau_0)
\]
whenever the function $\tau\to\gamma(\omega;x,t;\tau)$ does not have a jump at $\tau_0$. Next, one shows that at fixed $\tau_0\geq 0$, the set of points $(x,t)\in \R^2$ such that $\gamma(\omega;x,t;\tau)$ has a jump at $\tau_0$, has Lebesgue measure zero. Specializing to $\tau_0=0$, one can show in addition that the set of all $(\omega,x,t)$ such that
$\gamma(\omega;x,t;\tau)$ has a jump at $\tau = 0$ is jointly measurable, and it has measure zero (as a subset of $\Omega\times\R^2$).

Taking these observations together one obtains the desired result:

\emph{Convergence of characteristic curves.
For $\P$-almost all $\omega\in\Omega$ and Lebesgue almost all $(x,t)\in \R^2$,}
\begin{equation}\label{eq:convchar1}
\lim_{N\to\infty}\gamma^{(N)}(\omega;x,t;0)= \gamma(\omega;x,t;0).
\end{equation}

Details of the proof can be found in (Baumgartner et al. 2017).

\subsection{Convergence of approximate solutions}
\label{subsec:convsolu}

At level $N$, the propagation velocity of the discrete stochastic Goupillaud medium is given by the random field $c^{(N)}(\omega;x)$, where
\begin{equation*} 
c^{(N)}(\omega;x) = \Delta x_k^{(N)}(\omega)/\Delta t^{(N)}
\end{equation*}
in each random layer $x_{k-1}^{(N)}(\omega) \leq x < x_{k}^{(N)}(\omega)$ with $x_{k}^{(N)}(\omega) = L(\omega; t_k^{(N)})$ as defined in Subsection\;\ref{subsec:stochmod}.
At level $N$, the transport equation (Equation\;\ref{eq:trans}) reads
\[
\p_t u^{(N)}(\omega;t,x) + c^{(N)}(\omega;x) \p_x u^{(N)}(\omega;t,x)= 0 
\]
with initial data
$
u^{(N)}(\omega;t,x) = u_0(x).
$
If $u_0$ is a continuously differentiable function (or more generally, a locally integrable function with locally integrable weak first derivative), then
\[
    u^{(N)}(\omega;t,x) = u_0\big(\gamma^{(N)}(\omega;x,t;0)\big)
\]
is a weak solution at each fixed $\omega$. Indeed, at fixed $\omega$, the transport coefficient $c^{(N)}(\omega;x)$ is a piecewise constant function, and the characteristic curves $\gamma^{(N)}(\omega;x,t;\tau)$ are piecewise linear, continuous functions. Thus
$u^{(N)}(\omega;\cdot,\cdot)$ is a continuous function, and it satisfies the transport equation in each layer.
Define
\[
   u(\omega;t,x) = u_0\big(\gamma(\omega;x,t;0)\big).
\]
The results from Subsection \ref{subsec:convchar} lead to the convergence of $u^{(N)}$ to $u$ in the following sense:

\emph{Convergence of approximate solutions. Let $u_0$ be as described above. Then
\[
\lim_{N\to\infty} u^{(N)}(\omega;t,x) = u(\omega;t,x)
\]
for $\P$-almost all $\omega$ and Lebesgue almost all $(x,t)$, and
\[
   \lim_{N\to\infty} {\EE}^\P\|u^{(N)} - u\|_{L^p(K)} = 0
\]
whenever $K$ is a compact subset of $\R^2$ and $1\leq p < \infty$. Here}
\[
  \|u\|_{L^p(K)} = \left(\iint_K |u(t,x)|^p(d t,d x)\right)^{1/p}.
\]

Indeed, the almost sure convergence follows from the continuity of $u_0$ and Equation\;\ref{eq:convchar1}. From there, the convergence in mean of the local $p$th integrals follows from Lebesgue's dominated convergence theorem.

Note that a priori there is no meaning for $u$ to be a solution of the transport equation (Equation\;\ref{eq:trans}) other than being a limit of approximate solutions.

\begin{figure}[htb]
\centering
\includegraphics[width = 0.9\linewidth]{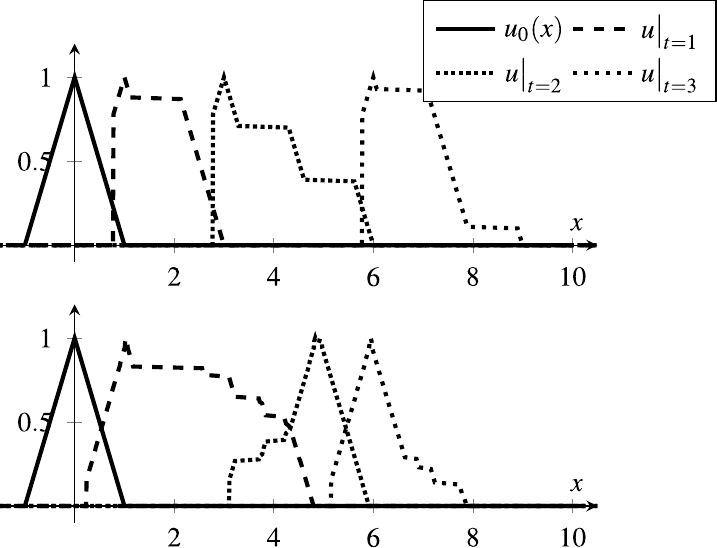}
		\caption{Figure 2. Trajectories of the  solution $u(x,t)$ at several time points. Top: The generating L\'evy process $L$ is Gamma distributed with scale parameter $k=1$, shape parameter $\theta=1$ and drift $d=1$. Bottom: The generating L\'evy process $L$ is Poisson distributed with jump size 1, intensity $c=1$ and drift $d=1$.}
\label{fig:U_realization}
\end{figure}

For the sake of illustration, the somewhat rough behavior of realizations of the limiting solutions are shown in Figure \ref{fig:U_realization}. The initial value $u_0$  is taken as a triangular function, the realizations of $u$ are shown at times $t = 1,2,3$. Two different L\'evy processes as drivers $L$ (cf. Subsection~\ref{subsec:stochmod}) are used, namely a Gamma process and a Poisson process, both with positive drift. The solutions have constant parts, which are created whenever the L\'evy process jumps at this point.

\section{Distribution of the characteristic curves}

Motivated by the results in the previous sections one may ask if it is possible to compute the statistical properties of the limit of the characteristic curves. Recall from Equations\;\ref{eq:charcurves} and \ref{eq:Lstar} that
\begin{equation}\label{eq:gammaxt0}
\gamma(x,t;0)=L(L^*(x)-t),
\end{equation}
where $L$ is a strictly increasing L\'evy process on a probability space $(\Omega,\nF,\PP)$ and $L^*(x)=\inf\menge{t\in\RR, L(t)\geq x}$. (To shorten notation, the random element $\omega$ will be omitted in the following.)
Since $\gamma(x,t;\tau)$ is just a translate of $\gamma(x,t;0)$, it suffices to compute the probability distribution of the latter.

For simplicity, only the case $t>0$ will be presented; for negative times one can do similar computations. In Subsection\;\ref{subsec:IG}, the explicit probability distribution of $L(L^*(x)-t)$ for the inverse Gaussian process will be elaborated. All the results are part of the work of the second author (Schwarz 2019), where detailed proof of the results of this section can be found.

\subsection{Results valid for any L\'evy process}
\label{subsection:any}

The considerations in this subsection are valid for any strictly increasing two-sided L\'evy process $L$ such that $L(0)=0$ almost surely.
For convenience of presentation, the notation of \emph{image measures} will be used in this section. The image measure of $\PP$ under a random variable $X$ will be denoted by $\bP_X$. In case of the random variable $L(L^*(x)-t)$  it reads
\[
\PP(L(L^*(x)-t)\in A) = \bP_{L(L^*(x)-t)} \ (A)
\]
where $A$ is any Borel measurable subset of $\RR$.
The distribution of ${L(L^*(x)-t)}$ can be explicitly computed if one knows when the L\'evy process hits $x$ (i.e., the \emph{hitting time} $L^*(x)$) and where it was right before (i.e., the \emph{undershoot} $L(L^*(x)^-) =\lim_{t\uparrow L^*(x)} L(t)$). By the law of total probability, the following decomposition holds:
\begin{align*}
	&\bP_{L(L^*(x)-t)} \ (A) \\[1pt]
	&= \iint_{\RR^2} \bP _{L(L^*(x)-t)\big| L^*(x)=s ,\hspace{1pt} L(L^*(x)^-)=y}(A)\bP_ {L^*(x),\hspace{1pt} L(L^*(x)^-)} (d s,d y)\\
	&= \iint_{\RR^2} \bP _{L(s-t)\big| L^*(x)=s ,\hspace{1pt} L(L^*(x)^-)=y}  (A) \bP_ {L^*(x),\hspace{1pt} L(L^*(x)^-)} (d s,d y),
\end{align*}
where
\begin{align*}
	&\bP _{L(s-t)\big| L^*(x)=s ,\hspace{1pt} L(L^*(x)^-)=y}  (A) = \PP(L(s-t)\in A| L^*(x)=s ,\hspace{1pt} L(L^*(x)^-)=y).
\end{align*}
So in order to obtain the distribution of $L(L^*(x)-t)$ one has to compute the joint probability distribution of the {hitting time} $L^*(x)$, the {undershoot} $L(L^*(x)^-)$, as well as the distribution of $L(s-t)$ for a given undershoot and hitting time. Detailed formulas for all possible locations of $x,y,s,t$ are collected in the following summary.

\emph{Summary of general formulas. Let $L$ be a strictly increasing L\'evy process on $\RR$  and $L^*(x)  =\inf\menge{t\in\RR: L(t)\geq x}$ be its hitting time. Let
	$ \bP_{L^*(x),\hspace{1pt} L(L^*(x)^-)} $
be the image measure of the joint probability distribution of the hitting time and the undershoot. Furthermore, let $t>0$ and let $A$ be a Borel measurable subset of $\RR$.}

\noindent
		\emph{(1) For $x\in\RR$ it generally holds that}
		\begin{align*}
			&\P(L(L^*(x)-t)\in A)\\
			&=\iint_{\RR^2} \bP _{L(s-t)\big| L^*(x)=s ,\hspace{1pt} L(L^*(x)^-)=y}  (A) \bP_ {L^*(x),\hspace{1pt} L(L^*(x)^-)} (d s,d y).
		\end{align*}
\noindent		
		\emph{(2) If $x>0$, $t\leq s$ and $0\leq y \leq x$, then
		\begin{align*}
			\bP_{L(s-t)\big| L^*(x)=s,\hspace{1pt} L(L^*(x)^-)=y}=\bP_{L'(s-t)\big| L'(s)=y}
		\end{align*}
		is the distribution of the so called \emph{L\'evy bridge}, constructed by means of an independent copy $L'$ from $L$,
        that is, $L'(r)| L'(s)=y$ is the process which starts at $0$ and reaches $y$ at time $s$, evaluated at time $r$.}
        %

\noindent		
		\emph{(3) For $x\geq0$, $0\leq s <t$ and $0\leq y \leq x$ one has
		$$\bP_{L(s-t)\big| L^*(x)=s,\hspace{1pt} L(L^*(x)^-)=y}(A) = \bP_{-L'(t-s)}(A),$$
		which is the probability distribution of an independent copy $L'$ of $L$.}

\noindent		
		\emph{(4) For $x<0$, $s<0$ and $y\leq x$ it holds that
		$$\bP_{L(s-t)\big| L^*(x)=s,\hspace{1pt} L(L^*(x)^-)=y}(A)= \bP_{y-L'(t)}(A),$$
		which is the probability distribution of $(y-L'(t))$, where $L'$ is an independent copy of $L$.}

Here are some indications how one may arrive at these results. Item (1) was deduced from the law of total probability above.
	
	The main difficulty in proving Item (2) is to show that
	\[ \bP_{L(r)\big| L(s) = y}=\bP_{L(r)\big| L(s)^- = y} \]
	for $r<s$, which means that the probability distribution does not change whether the undershoot $L(s)^-$ or the actual value of $L(s)$ is given. This can be proven by using the fact that a L\'evy process is, by definition, continuous in distribution.
	
	In Case (3) one has that $t>s$ and thus $s-t<0$, which results in a L\'evy process going into negative time direction and starting at zero.
	
\begin{figure}[htb] 	
\centering
	\includegraphics[width = 0.8\linewidth]{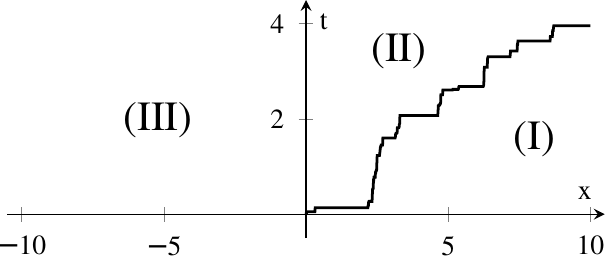}
		\caption{Figure 3. Regions corresponding to Cases (2)--(4). The plotted line is a trajectory of the L\'evy process $L$. In region (I), $x>0$ and $L^*(x)>t$ and hence this is Case (2). In region (II) $L^*(x)<t, x>0$, which corresponds to Case (3). Region (III) corresponds to Case (4).}
\label{Figure:Areas}
\end{figure}

	Case (4) is the same, but with a L\'evy process starting at the undershoot point $y$ and moving into negative time direction. The regions addressed in Cases (2)--(4) are visualized in Figure \ref{Figure:Areas}.

\subsection{The inverse Gaussian process}
\label{subsec:IG}

\begin{figure}[htb] 	
\centering
	\includegraphics[width = 0.9\linewidth]{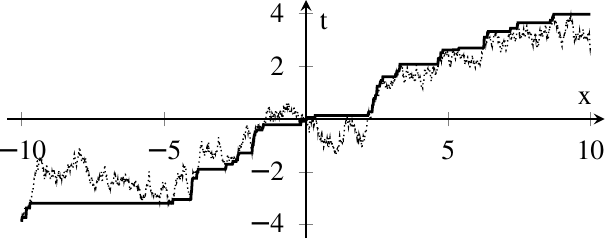}
		\caption{Figure 4. Illustration of running maxima. Dashed line: Brownian motion. Solid line: The two sided running maximum $M(x)$.}
\label{fig:construction_inverse_gaussian_process}
\end{figure}

This subsection is devoted to computing the explicit distribution of the characteristic curves in case of the inverse Gaussian process. The inverse Gaussian process can be constructed as follows: Let $(W^1(x))_{x\in [0,\infty)}$, $ (W^2(x))_{x\in [0,\infty)}$ be standard Brownian motions in space direction. \emph{Running maxima} are defined by
\begin{align*}
M^1(x)&=\max\menge{W^1(y),y\in [0, x]}\\
M^2(x)&=\min\menge{-W^2(y),y \in [0,-x ]}
\end{align*}
The \emph{two sided running maximum} is given by
\begin{align*}
	M(x)=\left\lbrace\begin{aligned}
		&M^1(x)&& \text{if } x\geq0\\
		&M^2(x)&& \text{if } x<0.
	\end{aligned}\right.
\end{align*}
An illustration can be found in Figure \ref{fig:construction_inverse_gaussian_process}.
Then the \emph{inverse Gaussian process} is defined by the hitting time of $M$:
\[ I(t)=\inf\menge{x\in\RR: M(x)> t}.\]
The probability density function of the inverse Gaussian process for $t\neq0$ can be found, e.g., in (Seshadri 1999) and is given by
\begin{align} \label{eqn:density_of_IG}
	f_{I(t)}(x)=\begin{cases}
	 g(x,t)	& \text{if} \quad t>0\\
	 h(x,t) & \text{if} \quad t<0\\
	\end{cases}
\end{align}
where
\begin{align*}
	g(x,t) &= \1_{(x>0)}\ \frac{t}{2\sqrt{\pi}} x^{-3/2}\re^{-t^2/(4x)},\\
	h(x,t) &= \1_{(x<0)} \frac{\abs{t}}{2\sqrt{\pi}} \abs{x}^{-3/2}\re^{-t^2/(4\abs{x})}.
\end{align*}
Here and in the sequel $\1_{A}$ denotes the characteristic function of an event $A$, that is, $\1_{A} = 1$ if $A$ occurs and $\1_{A} = 0$ otherwise.
In the following, the probability distribution components listed in the summary in Subsection \ref{subsection:any} will be computed. Due to the fact that the inverse Gaussian process has a probability density, all the components will be given by probability density functions, too.

Note that the distribution of $I(I^*(x)^-)$ is the same as $-I(I^*(-x))$, which means that the undershoot of $x$ is the negative \emph{overshoot} of $-x$. So, in order to compute the undershoot for $x<0$ is suffices to know the distribution of the overshoot for $-x$. Three preliminary steps are required.

\emph{The joint probability density of the hitting time, the undershoot, and the overshoot. Let $x\geq 0$ and $I$ be an inverse Gaussian process for positive times and let $I^*(x)$ denote the hitting time of the process. The joint probability density of the hitting time, the undershoot, and the overshoot is}
	\begin{align*}
		&f_{I^*(x),\hspace{1pt} I(I^*(x)^-),\hspace{1pt} I(I^*(x))}( s, a, b) \nonumber \\
		& = \1_{ (s\geq0, 0\leq a\leq x \leq b)} \frac{s}{2 \pi\sqrt{ a^3 (b-a)^3}} \exp\left(-\frac{s^2}{2a}\right).		
	\end{align*}

Indeed, it follows from the definition of the inverse Gaussian process that	
	\begin{align*}
		&I^*(x)=M(x),\\
		&I(I^*(x)^-)= \inf\menge{a \leq x: W(a)=M(x)}, \\
		&I(I^*(x))=\inf\menge{b\geq x: W(b) > M(x) }.
	\end{align*}
The following explicit formula for the joint probability density of $W$, $M$, and the undershoot is taken from (Karatzas 1988, Section 2.8, Proposition 8.15 and Remark 8.16):
	\begin{align*} 
		&f_{W(x),\hspace{1pt} M(x),\hspace{1pt} I(I^*(x)^-)}(  r,  s,  a )\\
		&=\1_{(s\geq r, s\geq 0, x\geq a \geq 0 )}\frac{s(s-r)}{\pi \sqrt{a^3(x-a)^3}} \exp\left( -\frac{s^2}{2a}-\frac{(s-r)^2}{2(x-a)}\right).
	\end{align*}
Using the reflection principle (Karatzas 1988, Section 2.6 A) this can be extended to
	\begin{align*}
	&f_{W(x),\hspace{1pt} M(x),\hspace{1pt} I(I^*(x)^-),\hspace{1pt} I(I^*(x))}(  r,  s,  a ,b)\\
	& = f_{W(x),\hspace{1pt} M(x),\hspace{1pt} I(I^*(x)^-)}(  r,  s,  a )\,\1_{(b\geq x)} \frac{s-r}{\sqrt{2 \pi (b-x)^3}}\exp\left(-\frac{(s-r)^2}{2(b-x)}\right).
	\end{align*}
	Finally, integrating over $r$ proves the desired result.

This result can now be used for providing the joint probability density of the hitting time and the undershoot.

\emph{The joint probability density of the undershoot and the overshoot.
	Let $I$ be a two-sided inverse Gaussian process. Then:}

\noindent
	\emph{(1) For $x>0$, $0\leq s$ and $0\leq y\leq x$ one has}
	\begin{align} \label{eq:Gauss1}
		&f_{I^*(x),\hspace{1pt} I(I^*(x)^-)} (s, y)  = \frac{s}{\pi\sqrt{y^3 (x-y)}} \exp\bigg(-\frac{s^2}{2y}\bigg).
	\end{align}
\emph{(2) For $x<0$, $s<0$ and $y\leq x$ one has}
	\begin{align} \label{eq:Gauss2}
			&f_{I^*(x),\hspace{1pt} I(I^*(x)^-)} ( s, y) =\int_x^0 \frac{-s}{2 \pi  \sqrt{a^3 (y-a)^3}} \exp\bigg(\frac{s^2}{2a}\bigg)  d a.
	\end{align}
\emph{(3) In all other cases $f_{I^*(x),\hspace{1pt} I(I^*(x)^-)} ( s, y) =0$.}

Indeed, for $x>0$ one starts from the formula
	\begin{align*}
		&f_{I^*(x),\hspace{1pt} I(I^*(x)^-)} ( s, y) = \int_x^\infty f_{I^*(x),\hspace{1pt} I(I^*(x)^-),\hspace{1pt} I(I^*(x))}( s, y, b) d b.
	\end{align*}
For $x<0$ one may use the fact that
\begin{align*}
&f_{I^*(x),\hspace{1pt} I(I^*(x)^-)} ( s, y) = f_{I^*(-x),\hspace{1pt} I(I^*(-x))} ( s, -y),
\end{align*}
which has been noted before. Again, integrating out the overshoot one obtains the desired results.

The next task is to compute the probability density of the inverse Gaussian bridge at time $r$, which starts at $0$ and reaches $y$ at time $s$.

\emph{The probability density of the inverse Gaussian bridge.
	Let $I$ be an inverse Gaussian process. Furthermore, let $s>r>0$, $y>0$ and let $f_{I(r)}(z)$ be the probability density of $I(r)$. Then}
	\begin{equation}\label{eq:Gauss3}
f_{I(r)\big| I(s)=y}(z)= \frac{f_{I(r)}(z) f_{I(s-r)} (y-z)}{f_{I(s)}(y)}.
    \end{equation}
	
Indeed, using the definition of conditional probabilities one has
that $f_{I(r)\big| I(s)=y}(z)$ equals
\begin{align*}
\frac{f_{I(r),\hspace{1pt} I(s)}(z,y)}{f_{I(s)}(y)}=\frac{f_{I(r),\hspace{1pt} I(s)-I(r)}(z,y-z)}{f_{I(s)}(y)}.
\end{align*}
Observing that $I(s)-I(r)$ is independent of $I(r)$, the term above is seen to be equal to
\begin{align*}
	\frac{f_{I(r)}(z) f_{I(s)-I(r)}(y-z)}{f_{I(s)}(y)}
			=\frac{f_{I(r)}(z) f_{I(s-r)}(y-z)}{f_{I(s)}(y)}
\end{align*}
where the last equality holds because $I$ is a L\'evy process and so $I(s)-I(r)$ has the same probability distribution as $I(s-r)$.

These results can now be collected, yielding the main result of this subsection. Recall from Equation\;\ref{eq:gammaxt0} that $\gamma(x,t;0)$ is given by
$I(I^*(x)-t)$ when $L = I$ is an inverse Gaussian process.

\emph{The probability density of the characteristic curve through the origin, inverse Gaussian process case.	
	Let $I$ be an inverse Gaussian process,
	\[f_{I(t)}(z)=\1_{(z\geq 0)} \left(\frac{t}{2\sqrt{\pi}}\right)z^{-3/2}\re^{-t^2/(2z)}  \]
and let $I'$ be an independent copy of $I$ and $t>0$. Then}

\noindent	
\emph{(1) the probability density $f_{I(I^*(x)-t)}(z)$ is given by}
	\begin{align*}
		&f_{I(I^*(x)-t)}(z)\\
     &=\int_0^\infty\int_{-\infty}^\infty f_{I(s-t) \big| I^*(x)=s,\hspace{1pt} I(I^*(x)^-)=y}(z)  f_{I^*(x),\hspace{1pt} I(I^*(x)^-)}(s,y)  d y d s,
	\end{align*}
\emph{where $f_{I^*(x),\hspace{1pt} I(I^*(x)^-)}(s,y)$ is given in Equations\;\ref{eq:Gauss1} and \ref{eq:Gauss2} and
	\[
    f_{I(s-t) \big| I^*(x)=s,\hspace{1pt} I(I^*(x)^-)=y}(z) = f_{I'(s-t)\big | I'(s)=y}(z)
    \]
is given in Equations\;\ref{eq:Gauss3} and \ref{eqn:density_of_IG}.}

\noindent	
\emph{(2) For $x\geq y>0$ and $s>t\geq0$,
	\begin{align*}
		&f_{I(s-t)\big | I^*(x)=s ,\hspace{1pt} I(I^*(x)^-)=y}(z)\\
		& = f_{I(s-t)\big | I(s)=y}(z)=\frac{f_{I(s-t)}(z) f_{I(t)} (y-z)}{f_{I(s)}(y)}.
	\end{align*}}
\noindent	
\emph{(3) For $x\geq0$ and $t > s>0$,
	\begin{align*}
		&f_{I(s-t)\big| I^*(x)=s ,\hspace{1pt} I(I^*(x)^-)=y}(A) = f_{I(t-s)}(-z).
	\end{align*}}
\noindent	
\emph{(4) For $x< 0$ and $s<0$,
	\begin{align*}
		&f_{I(s-t)\big| I^*(x)=s ,\hspace{1pt} I(I^*(x)^-)=y}(z) = f_{I(t)}(-y-z).
	\end{align*}}

\begin{figure*}[t]
	\includegraphics[width=\linewidth]{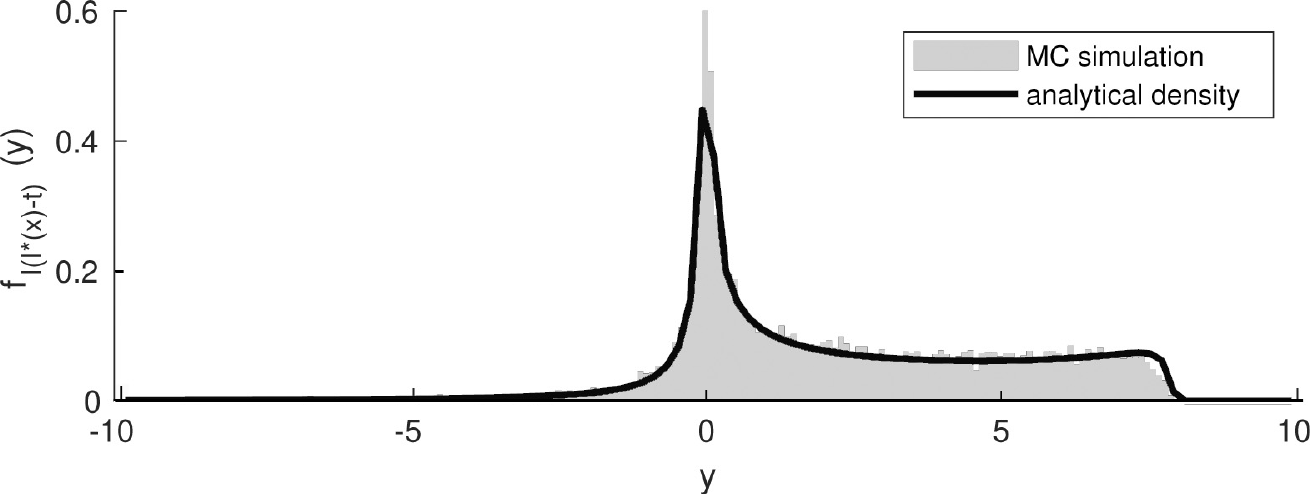}
	\caption{Figure 5. Empirical and analytical probability density of the base points of $\gamma(8,1;0) = I(I^*(8)-1)$.}
    \label{fig:iisx-t}
\end{figure*}

\emph{Example.} Recall that $\gamma(x_0,t_0;0) = I(I^*(x_0)-t_0)$ is the base point (intersection with the $x$-axis) of the characteristic curve that passes through $x_0$ at time $t_0$. Figure \ref{fig:iisx-t} shows the probability density $f_{I(I^*(x_0)-t_0)}(y)$ of these base points when $x_0=8$ and $t_0=1$. The black curve is the analytical probability density computed by the presented method. The grey bars represent a histogram of a Monte Carlo simulation of $1000$ inverse Gaussian processes and evaluated at $I^*(x_0)-t_0$. One sees that the probability density is zero for $y>8$ and that there is a concentration at $y=0$.

\section{Conclusion}

A Goupillaud medium is a piecewise constant layered medium such that the thickness of each layer is proportional to the corresponding propagation speed. A set-up has been developed for a specific stochastic Goupillaud medium in which the propagation speeds (or equivalently the thickness of the layers) are given by infinitely divisible random variables. Using a dyadic refinement, these random variables could be constructed as increments of a strictly increasing L\'evy process. It was shown that the one-dimensional transport equation can be solved in such a medium, and that the characteristic curves converge to shifted trajectories of the underlying L\'evy process as the time step goes to zero. If the initial data are sufficiently regular, the corresponding solutions converge pathwise and in the $p$th mean to a limiting function. It was shown how the probability distribution of the characteristic curves can be computed, explicitly calculated in the inverse Gaussian case. Modelling transport in discontinuities media by introducing random layers has a long tradition in applications. The paper shows that, through the proposed limiting procedure, transport in random media can be analyzed in media with a much higher degree of irregularity than previously available in the literature.

\section*{References}

\noindent
\setlength{\hangindent}{2em}
{\ten
Baumgartner, F., Oberguggenberger, M. \& Schwarz, M. 2017.
\newblock Transport in a stochastic {G}oupillaud medium.
\newblock In M.~Oberguggenberger, J.~Toft, J.~Vindas \& P.~Wahlberg (Eds.),
{Generalized Functions and Fourier Analysis}, Cham: Springer
  International Publishing, pp. 19-30.}

\noindent
\setlength{\hangindent}{2em}
{\ten
Burridge, R., Papanicolaou, G.S. \& White, B.S. 1998.
\newblock {One-dimensional wave propagation in a highly discontinuous medium.}
\newblock {International Journal of Approximate Reasoning}, 43:241-267.}

\noindent
\setlength{\hangindent}{2em}
{\ten
Flandoli, F. 2011.
\newblock {Random Perturbation of {PDE}s and Fluid Dynamic Models}. {Lecture Notes in Mathematics, Vol. 2015}.
\newblock Heidelberg: Springer.}

\noindent
\setlength{\hangindent}{2em}
{\ten
Fouque, J.P., Garnier, J., Papanicolaou, G. \& S{\o}lna, K. 2007.
\newblock {Wave propagation and time reversal in randomly layered media}. {Stochastic Modelling and Applied Probability, Vol. 56}.
\newblock New York: Springer.}

\noindent
\setlength{\hangindent}{2em}
{\ten
Ghanem, R.G. \& and Spanos, P.D. 1991.
\newblock {Stochastic Finite Elements: a Spectral Approach}.
\newblock New York: Springer-Verlag.}

\noindent
\setlength{\hangindent}{2em}
{\ten
Goupillaud, P.L. 1961.
\newblock An approach to inverse filtering of near-surface layer effects from
  seismic records.
\newblock {Geophysics}, 26:754-760.}

\noindent
\setlength{\hangindent}{2em}
{\ten
Hairer, M. 2014.
\newblock A theory of regularity structures.
\newblock Inventiones Mathematicae 198(2):269-504.}

\noindent
\setlength{\hangindent}{2em}
{\ten
Karatzas, I. \& Shreve, S.E. 1988.
\newblock Brownian Motion and Stochastic Calculus.
\newblock New York: Springer.}

\noindent
\setlength{\hangindent}{2em}
{\ten
Matthies, H.G. 2008.
\newblock Stochastic finite elements: computational approaches to stochastic
  partial differential equations.
\newblock Journal of Applied Mathematics and Mechanics (ZAMM) 88(11):849-873.}

\noindent
\setlength{\hangindent}{2em}
{\ten
Nair, B. \& White, B.S. 1991.
\newblock High-frequency wave propagation in random media---a unified approach.
\newblock SIAM Journal on Applied Mathematics 51(2):374-411.}

\noindent
\setlength{\hangindent}{2em}
{\ten
Sato, K. 1999.
\newblock {L}{\'e}vy Processes and Infinitely Divisible Distributions.
\newblock Cambridge: Cambridge University Press.}

\noindent
\setlength{\hangindent}{2em}
{\ten
Schwarz, M. 2019.
Stochastic Fourier Integral Operators and	Hyperbolic Differential Equations in Random Media; PhD thesis. Innsbruck: University of Innsbruck.}

\noindent
\setlength{\hangindent}{2em}
{\ten	
Seshadri., V. 1999.
The Inverse Gaussian Distribution: Statistical Theory and Applications. New York: Springer.}

\noindent
\setlength{\hangindent}{2em}
{\ten
Velo, A.P. \& Gazonas, G.A. 2019.
\newblock Applications of $z$-transforms to impact problems in layered elastic
  media.
\newblock {Archive of Applied Mechanics}, 89:581-590.}

\end{document}